\documentclass[12pt]{article}
\usepackage{latexsym, amssymb, amscd, amsfonts, epsfig, graphicx, colordvi,verbatim,ifpdf}
\usepackage{amsmath}
\usepackage{epsfig,cite, psfrag,eepic,color}
\usepackage{amscd,graphics}
\usepackage{graphicx}
\usepackage{fancybox}
\usepackage{float}
\usepackage{pstricks}
\usepackage{mathrsfs}
\usepackage{cases}
\usepackage[norsk,english]{babel}

\newtheorem{thm}{Theorem}[section]

\newtheorem{lem}[thm]{Lemma}

\def\pf{\noindent{\it Proof.} }
\def\qed{\nopagebreak\hfill{\rule{4pt}{7pt}}
\medbreak}

\numberwithin{equation}{section}

\def\qed{\nopagebreak\hfill{\rule{4pt}{7pt}}
\medbreak}

\newlength{\boxedparwidth}
  {\begin{center} \begin{tabular}{|@{\hspace{.315in}}c@{\hspace{.15in}}|}
                  \hline \\ \begin{minipage}[t]{\boxedparwidth}
                  \setlength{\parindent}{.25in}}%
  {\end{minipage} \\ \\ \hline \end{tabular} \end{center}}

\begin{document}
\begin{center}
{\large \bf Congruences on the Number of Restricted $m$-ary Partitions}
\end{center}

\begin{center}
{Qing-Hu Hou}$^{a}$, {Hai-Tao Jin}$^b$, {Yan-Ping Mu}$^{c}$ and {Li Zhang}$^{d}$ \vskip 2mm

   $^{a}$Center for Applied Mathematics\\
Tianjin University,  Tianjin 300072, P. R. China\\[6pt]
$^b$School of Science,
   Tianjin University of Technology and Education \\
   Tianjin 300222, P. R. China \\[6pt]
   $^{c}$College of Science, Tianjin University of Technology\\
 Tianjin 300384, P. R. China\\[6pt]
   $^{d}$Center for Combinatorics, LPMC-TJKLC\\
   Nankai University, Tianjin 300071, P. R. China\\

\vskip 2mm

  Emails:  qh\_hou@tju.edu.cn, jinht1006@tute.edu.cn, yanping.mu@gmail.com,
    zhangli427@mail.nankai.edu.cn
\end{center}

\vskip 6mm \noindent {\bf Abstract.} Andrews, Brietzke, R\o dseth and Sellers proved an infinite  family of congruences on the number of the restricted $m$-ary partitions when $m$ is a prime. In this note, we show that these congruences hold for arbitrary positive integer $m$  and thus confirm the conjecture of Andrews, et al.



\section{Introduction}
Let $m \ge 2$ be a fixed integer. An $m$-ary partition of a nonnegative integer $n$ is a partition of $n$ such that each part is a power of $m$. If there is ``no gaps'' in the parts, i.e., whenever $m^i$ is a part, $1,m,m^2,\ldots,m^{i-1}$ are parts, then the partition is called a restricted $m$-ary partition. The number of restricted $m$-ary partitions of $n$ is denoted by $c_m(n)$. Notice that the generating function of $c_m(n)$ is given by
\[
C_m(q):=\sum_{n=0}^{\infty}c_m(n)q^n
= 1+\sum_{n=0}^{\infty} \frac{q^{1+m+\cdots+m^n}}{(1-q)(1-q^m)\cdots(1-q^{m^n})}.
\]

In recent years, the arithmetic properties of $m$-ary partitions and restricted $m$-ary partitions have been widely studied since Churchhouse \cite{Churchhouse-1969} initiated the study of 2-ary partitions in the late 1960s. For example, R\o dseth \cite{R-1970} extended Churchhouse's results to include $p$-ary partition functions $b_p(n)$, where $p$ is any prime.   Andrews \cite{Andrews-1971}, Gupta \cite{Gupta-1972} and R\o dseth and Sellers \cite{R-S-2001} studied further the congruences for $b_m(n)$, where $m\ge 2$ is any positive integer. And later, Andrews, Fraenkel and Sellers \cite{Andrews-Fraenkel-Sellers-2015-1,Andrews-Fraenkel-Sellers-2015-2} provided characterizations of the values $b_m(mn)$ and $c_m(mn)$ modulo $m$.  Andrews, Brietzke, R\o dseth and Sellers \cite{Andrews-B-R-S-2015} proved that for odd prime $m$,
\[
c_m(m^{j+2}n+m^{j+1}+\cdots +m^2)\equiv 0 \pmod{m^j},
\]
for all $n \ge 0$ and $0 \le j <m$.
In this note, we will show that
\begin{thm}\label{conj}
For a fixed integer $m\ge 2$ and for all nonnegative integer $n$, we have
\[
c_m(m^{j+2}n+m^{j+1}+\cdots +m^2)\equiv 0 \pmod{\frac{m^j}{c_j}},
\]
where $c_j=1$ if $m$ is odd and $c_j=2^{j-1}$ if $m$ is even.
\end{thm}
This confirms the conjecture of Andrews et al.  \cite[Conjecture 4.1]{Andrews-B-R-S-2015}.

\section{Proof of Theorem 1.1}
Our proof is based on the following results and notations of Andrews et al. \cite{Andrews-B-R-S-2015}. Suppose that
\[
 {km+1 \choose j} = \sum_{i=0}^j s_{j,i}{k\choose i},
\]
where $s_{j,i}$ is an integer independent of $k$ such that $s_{j,0}=0$ if $j>1$ and $s_{j,0}=1$ otherwise. And denote
\[
h_j = \frac{q^j}{(1-q)^{j+1}}.
\]
\begin{lem}[Lemma 3.9 of \cite{Andrews-B-R-S-2015}]
\label{A-cm-j} Let $\overline{C_m}(q):=1-C_m(q)$, then we have
\begin{equation}\label{eq-cm-j}
\sum_{n=0}^{\infty}c_m(m^{j+2}n+m^{j+1}+\cdots+m^2)q^n
= \sum_{i=0}^{j+1}(P_{j,i}-Q_{j,i})h_i + \overline{C_m}(q)\sum_{i=1}^{j+1}(R_{j,i}-T_{j,i})h_i,
\end{equation}
where, for $i>0$, $P_{j,i}$ and $R_{j,i}$ are sums of monomials in $s_{u,v}$ of degree $j+1$ while $Q_{j,i}$ and $T_{j,i}$ are of degree $j$. In addition, $Q_{j,j+1}=Q_{j,0}=0$ and $P_{j,0}$ is of degree $j$ in the $s_{u,v}$. Furthermore, for $t>0$, we have the following recurrence relations:
\begin{align}
  &P_{j+1,t}= \sum_{i=t}^{j+2}(P_{j,i}+R_{j,i-1})s_{i,t}, \qquad
   R_{j+1,t}= \sum_{i=t}^{j+2}R_{j,i-1}s_{i,t},  \label{PR} \\[5pt]
  &Q_{j+1,t}= \sum_{i=t}^{j+2}(Q_{j,i}+T_{j,i-1})s_{i,t}, \qquad
   T_{j+1,t}= \sum_{i=t}^{j+2}T_{j,i-1}s_{i,t}, \label{QT} \\[5pt]
  &P_{j+1,0}= P_{j,1}=Q_{j+1,1}. \label{P0}
\end{align}
\end{lem}

We firstly introduce an evaluation on monomials in $s_{u,v}$ by defining
\[
v(s_{j,j})= j, \quad v(s_{j,j-1}) = j-1-\varepsilon, \quad v(s_{j,i})=0,\ \forall\, 1\le i<j-1,
\]
and
\[
v(s_{j,i} s_{j',i'} \cdots) = v(s_{j,i}) + v(s_{j',i'}) + \cdots,
\]
where $\varepsilon$ is an arbitrarily small  positive real number.
\begin{lem}\label{div}
Let $p = s_{j,i} s_{j',i'} \cdots$ be a monomial in $s_{u,v}$. If
$v(p) > \ell$,  then $p$ is divisible by $\frac{m^{\ell+1}}{c_{\ell+2}}$.
\end{lem}
\pf
It is known that $s_{i,i}=m^i$. We will show that
\begin{equation}
s_{i+1,i}=\frac{i(m-1)}{2}m^i+m^i. \label{si+1i}
\end{equation}

Recall that
\begin{equation}\label{j>1-1}
  {km+1 \choose j} = \sum_{i=1}^j s_{j,i}{k\choose i}.
\end{equation}
Since
\begin{align}\label{lhs}
  {km+1 \choose j}
&= \frac{(km+1)km(km-1)\cdots(km-j+2)}{j!} \qquad \nonumber \\[5pt]
&=\frac{m^jk^j}{j!}-\frac{(j-3)}{2}\frac{m^{j-1}k^{j-1}}{(j-1)!}+\cdots,
\end{align}
and
\begin{align}\label{rhs}
  \qquad \sum_{i=1}^j s_{j,i}{k\choose i}
&= s_{j,j}{k\choose j} + s_{j,j-1}{k\choose j-1} + \cdots \nonumber \\[5pt]
&= s_{j,j}\frac{k^j}{j!}-s_{j,j}\frac{(j-1)}{2}\frac{k^{j-1}}{(j-1)!}
 + s_{j,j-1}\frac{k^{j-1}}{(j-1)!}+\cdots.
\end{align}
By  comparing the coefficients of $k^{j-1}$ in \eqref{lhs} and \eqref{rhs}, we get \eqref{si+1i} immediately.

Now suppose that
\[
p = s_{1,1}^{a_1} s_{2,2}^{a_2} \cdots s_{2,1}^{b_1}  s_{3,2}^{b_2} \cdots s_{3,1}^{c_1}  \cdots.
\]
Let $\ell = \sum_{i} i a_i + \sum_j  j b_j$.
We have
\[
\ell \ge  v(p) = \sum_{i} i a_i + \sum_j (j-\varepsilon) b_j >  \ell - 1.
\]
If $m$ is odd, $p$ is divisible by $m^\ell$. If $m$ is even, $p$ is divisible by
\[
m^{\sum_{i} i a_i} \cdot \left(\frac{m}{2} \right)^{b_1} \cdot \left(\frac{m^2}{2} \right)^{b_2} \cdots = \frac{m^\ell}{2^{\sum_j b_j}}.
\]
In both cases, we have $\frac{m^{\ell}}{c_{\ell+1}} \mid p$, which completes the proof. \qed

Now we consider the monomials in $P_{j,i}, Q_{j,i}, R_{j,i}$ and $T_{j,i}$. We will characterize the terms with minimal evaluations.

\begin{lem} There exists a unique term with minimal evaluation in $P_{j,i}$ ($Q_{j,i}, R_{j,i}$ and $T_{j,i}$, respectively). Denote such a term by $\overline{P_{j,i}}$ ($\overline{Q_{j,i}},  \overline{R_{j,i}}, \overline{T_{j,i}}$, respectively).
We have
\begin{itemize}
  \item[{\rm(1)}] For $1\le i \le j+1$,
  \begin{equation}\label{P-R}
   \overline{P_{j,i}}=\overline{R_{j,i}}, \quad
   v(\overline{P_{j,i+1}})-v(\overline{P_{j,i}})>i;
  \end{equation}
  \item[{\rm(2)}]For $1\le i \le j$,
  \begin{equation}\label{Q-T}
   \overline{Q_{j,i}}=\overline{T_{j,i}}, \quad
  v(\overline{Q_{j,i+1}})-v(\overline{Q_{j,i}})>i.
  \end{equation}
  \end{itemize}
\end{lem}
\pf We will show that \eqref{P-R} and \eqref{Q-T} hold by induction on $j$. The uniqueness of the minimal term (i.e., the term with minimal evaluation) follows from the proof.
\begin{itemize}
  \item[{\rm(1)}]
When $j=1$, by Lemma 3.7 of \cite{Andrews-B-R-S-2015}, we have
\[
 \overline{P_{1,1}} = \overline{R_{1,1}} = s_{1,1}s_{2,1}, \quad
 \overline{P_{1,2}} = \overline{R_{1,2}} = s_{1,1}s_{2,2}.
\]
Thus \eqref{P-R} holds by straightforward checking.

Now we assume that \eqref{P-R} holds for positive integer $j$ and we seek for the minimal term in $P_{j+1,i}$.
From the recurrence relation \eqref{PR}, we see that for $1\le i \le j+2$,
\begin{equation}\label{rec-p}
  P_{j+1,i}
  = P_{j,i}s_{i,i}+P_{j,i+1}s_{i+1,i}+\cdots
   + P_{j,j+1}s_{j+1,i} +R_{j,i-1}s_{i,i}+R_{j,i}s_{i+1,i}
   + \cdots + R_{j,j+1}s_{j+2,i},
\end{equation}
and
\begin{equation}\label{rec-r}
  R_{j+1,i}
  = R_{j,i-1}s_{i,i}+R_{j,i}s_{i+1,i}
   +\cdots + R_{j,j+1}s_{j+2,i}.
\end{equation}
By the induction hypothesis, we have that for $i+1\le k \le j+1$,
\begin{equation*}
v(\overline{P_{j,k}}) - v(\overline{P_{j,i}}s_{i,i})
= v(\overline{P_{j,k}}) - v(\overline{P_{j,i}})-i>0,
\end{equation*}
and
\[
v(\overline{R_{j,k}}) - v(\overline{R_{j,i}}s_{i+1,i})
= v(\overline{R_{j,k}}) - v(\overline{R_{j,i}})-i+\varepsilon>0.
\]
Therefore, the minimal term must lie in $P_{j,i}s_{i,i}, R_{j,i-1}s_{i,i}$ or $R_{j,i}s_{i+1,i}$.
If $i > 1$, then
\[
v(\overline{P_{j,i}}s_{i,i}) > v(\overline{R_{j,i}} s_{i+1,i}) > v(\overline{R_{j,i-1}} s_{i,i})  + i - 1  - \varepsilon > v(\overline{R_{j,i-1}} s_{i,i}).
\]
If $i=1$, then
\[
v(\overline{P_{j,i}}s_{i,i}) > v(\overline{R_{j,i}} s_{i+1,i}).
\]
Hence
\begin{equation}\label{rel-PR}
  \overline{P_{j+1,i}}=\overline{R_{j+1,i}}
=\left\{
       \begin{array}{ll}
        \overline{R_{j,i-1}}s_{i,i}, & \hbox{if}\quad 2\le i\le j+2, \\[7pt]
        \overline{R_{j,1}}s_{2,1}, & \hbox{if} \quad i=1.
        \end{array}
  \right.
\end{equation}
This implies that for $i>1$,
\begin{align*}
v(\overline{P_{j+1,i+1}})-v(\overline{P_{j+1,i}})
&=v(\overline{R_{j,i}}s_{i+1,i+1})-v(\overline{R_{j,i-1}}s_{i,i})\\[5pt]
&=v(\overline{R_{j,i}})-v(\overline{R_{j,i-1}})+1\\[5pt]
&>i-1+1\\[5pt]
&=i,
\end{align*}
and for $i=1$,
\begin{align*}
v(\overline{P_{j+1,2}})-v(\overline{P_{j+1,1}})
&=v(\overline{R_{j,1}}s_{2,2})-v(\overline{R_{j,1}}s_{2,1})\\[5pt]
&=v(s_{2,2})-v(s_{2,1})\\[5pt]
&=2-(1-\varepsilon)\\[5pt]
&>1.
\end{align*}
Thus \eqref{P-R} follows by mathematical induction.

\item[{\rm(2)}]
By a discussion similar to (1), we obtain \eqref{Q-T}. Moreover, we
have the following recurrence relation:
\begin{equation}\label{rel-QT}
  \overline{Q_{j+1,i}}=\overline{T_{j+1,i}}
=\left\{
       \begin{array}{ll}
        \overline{T_{j,i-1}}s_{i,i}, & \hbox{if}\quad 2\le i\le j+2, \\[7pt]
        \overline{T_{j,1}}s_{2,1}, & \hbox{if} \quad i=1.
        \end{array}
  \right.
\end{equation}
\end{itemize}
This completes the proof. \qed

Now we are ready to give a proof of Theorem \ref{conj}.

{\noindent \emph{Proof of Theorem \ref{conj}}}.
By iterative use of Equation \eqref{rel-PR}, we derive that for $1\le i\le j+1$,
\begin{equation}\label{res-pr}
 \overline{P_{j,i}}=\overline{R_{j,i}}
=\left\{
       \begin{array}{ll}
        s_{1,1}s_{2,2}\cdots s_{i,i}s_{2,1}^{j-i+1}, & \hbox{if}\quad 2\le i\le j+1, \\[7pt]
        s_{1,1}s_{2,1}^{j}, & \hbox{if} \quad i=1.
        \end{array}
  \right.
\end{equation}
Similarly, by Equation \eqref{rel-QT}, we deduce that for $1\le i\le j+1$,
\begin{equation}\label{res-qt}
  \overline{Q_{j,i}}=\overline{T_{j,i}}
=\left\{
       \begin{array}{ll}
        s_{1,1}s_{2,2}\cdots s_{i,i}s_{2,1}^{j-i}, & \hbox{if}\quad 2\le i\le j+1, \\[7pt]
        s_{1,1}s_{2,1}^{j-1}, & \hbox{if} \quad i=1.
        \end{array}
  \right.
\end{equation}
It is easy to see that $v(\overline{Q_{j,i}}) < v(\overline{P_{j,i}})$ and for $2\le i\le j+1$,
\begin{align*}
v(s_{1,1}s_{2,2}\cdots s_{i,i}s_{2,1}^{j-i})
- v(s_{1,1}s_{2,1}^{j-1})
= \frac{i(i-1)}{2} + \varepsilon(i-1) >0.
\end{align*}
Hence among all monomials in $P_{j,i}, Q_{j,i}, R_{j,i}$ and $T_{j,i}$, the minimal term is
$s_{1,1} s_{2,1}^{j-1}$. Let $p$ be a monomial in $P_{j,i}, Q_{j,i}, R_{j,i}$ or $T_{j,i}$. By induction we see that $s_{1,1}$ is a factor of $p$. Write $p=s_{1,1}p'$. Then
\[
v(p') \ge v(s_{2,1}^{j-1}) > j-2.
\]
By  Lemma~\ref{div} we derive that
\[
\frac{m^{j-1}}{c_j} \mid p',
\]
and hence $\frac{m^j}{c_j} \mid p$. Theorem~1.1 follows immediately from the generating function \eqref{eq-cm-j}. \qed

\vspace{.3cm}

{\noindent \bf Acknowledgments.}
 This work was supported by the 973 Project, the PCSIRT Project
of the Ministry of Education and the National Science Foundation of China.
The second author was also supported by the
National Science Foundation of China (No. 11501416).

\end{document}